\documentstyle[fleqn,twoside,plenum2,a_i_te2]{article}

\begin{document}

\def\uu{\accent23u}
\mathindent=2cm
\def\statement#1#2{\medskip{\bf #1}~~{\sl #2}\medskip}
\def\secondstatement#1#2{{\bf #1}~~{\sl #2}\medskip}
\def\I{\kern-1.5ex}
\def\D{\Delta}
\def\G{\Gamma}
\def\la{\lambda}
\def\e{\varepsilon}
\def\l{\ell}
\def\TT{\mathop{\cal T}\nolimits}
\def\FF{\mathop{\cal F}\nolimits}
\def\GG{\mathop{\cal G}\nolimits}
\def\eql    {\mathop{=}      \limits}
\def\cupl   {\mathop{\cup}   \limits}
\def\maxl   {\mathop{\max}   \limits}
\def\suml   {\mathop{\sum}   \limits}
\newsymbol\smallsetminus 2072
\def\setminus#1#2{#1\kern-.08ex\smallsetminus\kern-.08ex#2}
\def\_#1{\mathop{\hspace{-2pt}^{}_{#1}}}
\def\1n{1,\ldots,n}

\Btit{}

\tit{519.172}
    {P.~Yu.~Chebotarev and E.~V.~Shamis}
    {The Matrix-Forest Theorem and Measuring Relations in\\
     Small Social Groups}
    {\footnote[1]{This work was supported by the Russian Foundation
                  for Fundamental Research, Grant No.~96-01-01010.}}

\translation{Institute of Control Sciences, Russian Academy of
Sciences, Moscow. Translated from Avtomatika i Telemekhanika, No.~9,
pp.~125--137, September, 1997. Original article submitted December 23,
1996.}

\abst{We propose a family of graph structural indices related to the
matrix-forest theorem. The properties of the basic index that expresses
the mutual connectivity of two vertices are studied in detail. The
derivative indices that measure ``dissociation," ``solitariness," and
``provinciality" of vertices are also considered.  A nonstandard metric
on the set of vertices is introduced, which is determined by their
connectivity. The application of these indices in sociometry is
discussed.
}

\arcsect{1}{INTRODUCTION}

Given a graph, how should one evaluate the proximity between its
vertices? The standard distance function is the length of the shortest
path. But is it not worth taking into account the {\em number\/} of
paths between vertices? Which vertices can be considered central and
which peripheral? Which graphs are dense, and which are sparse? Which
are homogeneous? The choice of indices that express these and other
structural properties of graphs depends on applications, more exactly,
on the type of applications. This type should ideally be formulated in
terms of axiomatic requirements on the structural indices or via
modeling those concepts that should be evaluated by these indices. The
applications are numerous; essentially, these are all applications of
graph theory: transport, reliability, transmission of information,
structural modeling, chemistry, molecular biology, epidemiology, etc.

The application we shall focus on is one of the most difficult to
formalize. It is sociology, more precisely, sociometry where structural
indices are usually chosen heuristically.

Sociometry studies the structure of small social groups on the basis of
given relations on them. As a rule, these relations are binary; in some
cases they are weighted. Small social groups are groups where public
relations manifest themselves in the form of personal contacts or,
simply stated, these groups are natural communities where everyone
knows each other. The binary relations under study mainly result from
sociometric interrogations. For example, each member of a group is
asked to indicate those persons with whom she is in sympathy (or out of
sympathy), or with whom she spends her spare time most often, or with
whom she {\em would prefer\/} to cooperate in certain activities (work,
rest, ``exploration," etc.), or who, in her opinion, has certain
characteristics. If a member $i$ of a group indicates (among others)
$j$, an arc from $i$ to $j$ is drawn in the digraph of the
relationship.  Nonoriented graphs are frequently included to represent
objective information (contacts, collaborations, etc.). If a set of
similar questions is asked or the respondents report their assumptions
on the opinion of others (autosociometric data), multigraphs or
multidigraphs can serve as the model. A similar approach is used in
political studies where countries or parties involved in certain
relationships are investigated.

A lot of various kinds of relations can be studied, each requiring its
own properties of the structural indices, so it is problematic to
construct the desirable axiomatics for every application. Another
approach seems more realistic: to collect a ``library" of structural
indices with specified properties, and to use those indices whose
features are most appropriate for the relations under study.

Traditional sociological indices are very simple. For instance, the
{\em sociometric status\/} of the $i$th member of a group is the
normalized in-degree (the number of entering arcs) of the $i$th vertex;
the {\em psychological effusiveness\/} is the normalized out-degree;
the {\em reciprocity of the choice of $i$} is the normalized
number of pairs of opposite arcs incident to $i$ \cite{Pan}. The
{\em density\/} and the {\em cohesion\/} of a group result by averaging
the sociometric status and the reciprocity over the group. The {\em
heterogeneity\/} of a group is measured by the empirical variance of
the sociometric status over the group. The imperfection of these
elementary indices is caused by their local nature. In this connection,
Paniotto \cite{Pan} adduces examples of essentially dissimilar
structures with the same values of the above indices and states that
more sensitive indices that capture the topological structure as a
whole and the part of each individual in this structure are desirable.

A family of more sensitive indices is based on evaluating the
vertex status by the sum of the lengths (or reciprocal lengths) of the
shortest paths that lead to this vertex from all other vertices
\cite{HeGl}. Note, however, that such characteristics of a group member
are frequently unchanged on altering the connections between other
members, which may not conform with the interpretation of the model.
Furthermore, when characterizing the proximity of two members of a
group, it is often worth taking into account not only the length of the
shortest path between them, but also the numbers of paths of various
lengths.

One more idea employed in the construction of sociometric indices is
to evaluate the group cohesion by the number of arcs (edges) in the
minimum cutset, i.e., by the minimum number of connections whose
removal breaks the connectedness of the corresponding graph. The
normalized minimum number of members whose removal (together with their
connections) makes the graph disconnected is sometimes called
{\em group stability\/ {\rm (}vitality}). One problem (besides the
computational one) with such indices is that the graph can be found to
be disconnected from the very beginning. This situation still allows
one to study the increase in the number of connected components. On the
other hand, even for a connected graph, these indices are solely
determined by its ``bottlenecks," i.e., such characteristics are
indifferent to the existence of joined subgroups with relatively poor
connections between them.

In this paper, we propose a family of sensitive structural indices and
study its properties. The definitions of the indices are based on the
matrix-forest theorem (Section~2). Section~3 is devoted to the
properties of the basic index of vertex proximity; Section~4 discusses
it and introduces derivative indices.

The basic results of this work were stated in \cite{ChShAtl}. A
topological interpretation of the vertex proximity index (implicitly
used in \cite{Che89}) was obtained in \cite{ShaAMS,Sha}; an
interesting further investigation of the matrix of these indices was
undertaken in \cite{Merr}; close ideas with reference to chemistry were
developed in \cite{GoDrRo}, where an important analogy with electrical
networks was also formulated. In a subsequent paper, we are going to
compare the structural indices proposed here with other ones known from
the literature (see, for example,~\cite{Frie}).

\arcsect{2}{THE MATRIX-FOREST THEOREM}

The matrix-forest theorem is formulated for multigraphs and
multidigraphs (which differ from graphs and digraphs by the possibility
of multiple edges and arcs). A {\it subgraph\/} of a multigraph $G$ is a
multigraph all of whose vertices and edges belong to the vertex set and
the edge set of $G$. A {\it spanning subgraph\/} of a multigraph $G$ is a
subgraph of $G$ with the same vertex set as that of $G$. A {\it path\/}
in a multigraph $G$ is an alternating sequence of distinct vertices and
edges, which starts and ends with vertices and has each edge situated
between two vertices incident to it. Sometimes we consider a path
as a subgraph of $G$. A {\it forest\/} is a cycleless graph. A {\it
tree\/} is a connected forest. A {\it rooted tree\/} is a tree with one
marked vertex, called a {\it root}. Formally, a rooted tree is a pair
$(T,r)$, where $T$ is a tree and $r$ is its vertex. A {\it component\/}
of a multigraph $G$ is any maximal (by inclusion) connected subgraph of
$G$.  Obviously, all components of a forest are trees.

A {\it rooted forest\/} is defined as a forest with one marked vertex
in each component. A {\it directed path\/} in a multidigraph $\G$ is
defined similarly to a path in a multigraph, but here each arc is
directed from the previous vertex to the next one in the sequence. A
digraph is called a {\it directed tree\/} (a {\it directed forest}) if
the graph obtained from it by replacement of all its arcs with edges is
a tree (a forest). The definitions of {\it directed rooted tree\/} and
{\it directed rooted forest\/} are analogous to the definitions of
rooted tree and rooted forest (we will omit the word ``directed" while
talking about subgraphs of $\G$). A {\it diverging tree\/} is a
directed rooted tree that contains directed paths from the root to all
other vertices. A {\it diverging forest\/} is a directed rooted forest,
all of whose components are diverging trees.

Suppose that $G$ is a weighted multigraph with vertex set $V(G)=\{\1n\}$
and edge set $E(G)$. Let $\e_{ij}^p\ge0$ be the {\em weight\/} of the
$p$th edge between vertices $i$ and $j$ in $G$. This weight will be
also referred to as the {\it conductance\/} of the edge.

The {\em Kirchhoff matrix\/} of $G$ is
the $n\times n$ matrix $L=L(G)=\left(\l\_{ij}\right)$ with
\begin{eqnarray}
\label{lij}
\l\_{ij}\I&=&\I-\suml_{p=1}^{a\_{ij}}\e_{ij}^p,\;\:j\neq i,\;\:i,j=\1n,
\\
\l\_{ii}\I&=&\I-\suml_{j\neq i}\l\_{ij},\;\:i=\1n,
\label{lii}
\end{eqnarray}
where $a\_{ij}$ is the number of edges between $i$ and $j$. The product
of the weights of all edges that belong to a subgraph $H$ of a
multigraph $G$ will be referred to as the {\it weight\/} or {\it
conductance\/} of $H$ and denoted by $\e(H)$. The weight (conductance)
of a subgraph without edges is set to be~1. For every nonempty set of
subgraphs $\GG$, its weight is defined as follows:
\[
\e(\GG)=\suml_{H\in\GG}\e(H).
\]

The weight of the empty set is zero.

The following matrix-forest lemmas are similar to the classical
matrix-tree theorems, obtained by Kirchhoff and some other writers in
the nineteenth century (for the history, see \cite{Moo}). We shall
formulate Tutte's generalization of the matrix-tree theorem to weighted
multigraphs (see \cite{Tut}).

Denote by $L^{ij}$ the cofactor of $\l\_{ij}$ in $L$. Let $\TT(G)=\TT$
be the set of all spanning trees of multigraph $G$.

\statement{THEOREM~1~(matrix-tree theorem for weighted
multigraphs).}{For any weighted multigraph $G$ and for any $i,j\in
V(G),$ $L^{ij}=\e(\TT).$}

Tutte also obtained an analogous result for weighted multidigraphs.

Let $\G$ be a multidigraph with vertex set $V(\G)=\{\1n\}$, and suppose
that $\e_{ij}^p$ is the {\em weight\/} (or the {\em conductance\/}) of
the $p$th arc from $i$ to $j$ in $\G$. The Kirchhoff matrix of $\G$ is
the $n\times n$ matrix $L=L(\G)=\left(\l\_{ij}\right)$ with entries
$\l\_{ij}=-\suml_{p=1}^{a\_{ji}}\e_{ji}^p\:$, $j\neq i,\;\:i,j=\1n$,
and
$\l\_{ii}=-\suml_{j\neq i}\l\_{ij}\:$, $i=\1n$, where $a\_{ji}$ is the
number of arcs from $j$ to $i$ in $\G$. Observe that $\l\_{ii}$ is the
total conductance of the arcs {\it converging\/} to $i$. The
conductance (weight) of a subgraph of $\G$ and the weight of a set of
multidigraphs are defined analogously to the case of multigraphs.

Suppose that $\TT^i$ is the set of all spanning trees of $\G$ diverging
from $i$, and $L^{ij}$ is the cofactor of $\l\_{ij}$ in $L$, as before.

\statement{THEOREM~2~(matrix-tree theorem for weighted
multidigraphs).}{For any weighted multidigraph $\G$ and for any
$i,j\in V(\G),$ $L^{ij}=\e(\TT^i).$}

Observe that in the directed case, entries in different rows of $L$ may
have different cofactors, but all the entries of the same row have
equal cofactors. For simplicity, Tutte formulates these theorems only
for diagonal cofactors $L^{ii}$. The ``directed" matrix-tree theorem
concerning all $L^{ij}$ is given in \cite{HaPa}. If the weights of all
edges (arcs) are ones, Theorems~1 and~2 tell us about the {\it
numbers\/} of the corresponding spanning trees.

We shall now formulate the matrix-forest lemmas and the matrix-forest
theorem.

Consider the matrices
\[
W(G)=I+L(G)
\]
and
\[
W(\G)=I+L(\G),
\]
where $I$ is the identity matrix. $W^{ij}(G) $ and $W^{ij}(\G)$ will
denote the cofactors of the $(i,j)$-entries of $W(G)$ and $W(\G)$,
respectively.

Suppose that $\FF(G)=\FF$ is the set of all spanning rooted forests of a
weighted multigraph $G$ and $\FF^{ij}(G)= \FF^{ij}$ is the set of those
spanning rooted forests of $G$ such that $i$ and $j$ belong to the
same tree rooted at $i$. Let $W=W(G),\;W^{ij}=W^{ij}(G)$.

\statement{LEMMA~1~(matrix-forest lemma for weighted multigraphs).}
{For any weighted multigraph $G,$\\
{\rm (1)} $\det W=\e(\FF);$\\
{\rm (2)} for any $i,j\in V(G),$ $\;W^{ij}=\e(\FF^{ij}).$}

Since the matrix $W$ of a weighted multigraph is symmetric, item~(2) of
Lemma~1 remains true if we replace $\FF^{ij}$ by $\FF^{ji}$.

Suppose that $\FF(\G)=\FF$ is the set of all spanning diverging forests
of multidigraph $\G$ and $\FF^{i\rightarrow j}(\G)=\FF^{i\rightarrow
j}$ is the set of those spanning diverging forests of $\G$ such that
$i$ and $j$ belong to the same tree diverging from $i$. Let
$W=W(\G),\;W^{ij}=W^{ij}(\G)$.

\statement{LEMMA 2~(matrix-forest lemma for weighted multidigraphs).}
{For any weighted multidigraph~$\G,$\\
{\rm (1)} $\det W=\e(\FF);$\\
{\rm (2)} for any $i,j\in V(\G),$ $\;W^{ij}=\e(\FF^{i\rightarrow j}).$}

A similar lemma can be formulated for converging forests.

If the matrix $W^{-1}$ exists, we will denote it by
\begin{equation}
Q=(q\_{ij})=W^{-1}=(I+L)^{-1}
\label{DefQ}
\end{equation}
(either for a weighted multigraph $G$ or for a weighted multidigraph
$\G$). Then $Q=(\det W)^{-1}W^*$, where
$W^*=\left(W^{ij}\right)^{\scriptscriptstyle T}$ is the adjugate
matrix of $W$. The matrix-forest theorem \cite{ShaAMS,ChShAtl} follows
from Lemmas~1 and~2.

\statement{THEOREM~3~(matrix-forest theorem).}{\\
\indent {\rm 1.} For any weighted multigraph $G,$ the matrix $Q=W^{-1}$
exists and
$q\_{ij}=\e(\FF^{ji})\big/\e(\FF),\;\;i,j=\1n.$\\
\indent {\rm 2.} For any weighted multidigraph $\G,$ the matrix
$Q=W^{-1}$ exists and
$q\_{ij}=\e(\FF^{j\rightarrow i})\big/\e(\FF),\;\;i,j=\1n.$}

If the weights of all edges (arcs) are ones, the weights of sets of
spanning forests in Lemmas~1 and~2 and Theorem~3 are equal to the {\it
numbers\/} of the corresponding forests.

Lemma~2 can be derived in the shortest way from one version of
Chaiken's result \cite{Cha}, namely, by putting $U=W=\varnothing$ and
then $U=\{i\},\: W=\{j\}$ in the first formula on page~328 (cf.\
\cite[Theorem~3.1]{MoonDM}). A longer inference results by the
sequential application of results from \cite{Kel,KeCh,FiSe,MaOl}. This
also provides an interpretation for the inverse Laplacian
characteristic matrix of a multidigraph. An inference of Lemma~1 from
Lemma~2 is given in the Appendix, as well as the proofs of the
following results. Another complete (i.e., not exploiting any strong
theorems) proof of Lemma~1 for the case of equal weights of edges is
contained in \cite{Sha}.

The matrix-forest theorem allows us to consider the matrix $Q=W^{-1}$
as the matrix of {\em ``relative forest accessibilities"\/} (in short,
{\em accessibilities}) of the vertices of $G$ (or $\G$). These values
can be used to measure the proximity between vertices (the ``farther"
$i$ from $j$, the smaller is $q\_{ij}$). This interpretation is
validated by the properties presented in the following section. For
simplicity, these properties are formulated for nonoriented
multigraphs, although many of them have ``oriented" counterparts which
can be proved similarly.

\arcsect{3}{PROPERTIES OF THE RELATIVE FOREST ACCESSIBILITIES}
\label{SPrope}

Suppose that $G$ is a weighted multigraph with strictly positive
weights of edges, and let
\[
Q=(q\_{ij})=W^{-1}
\]
be its matrix of relative forest accessibilities.

\statement{PROPOSITION 1.}{For any $G,$ matrix $Q$ is symmetric.}

\secondstatement{PROPOSITION 2.}{For any $G,$ $Q$ is a doubly
stochastic matrix$,$ i.e.$,$
\newline{\rm (1)} $q\_{ij}\geq0, \;\: i,j=\1n;$
\newline{\rm (2)} $\suml_{j=1}^n q\_{ij}=1,\;\: i=\1n;$
\newline{\rm (3)} $\suml_{i=1}^n q\_{ij}=1,\;\: j=\1n.$
}

According to this property, $q\_{ij}$ may be interpreted as the
fraction of the connectivity of vertices $i$ and $j$ in the total
connectivity of $i$ with all vertices.

\statement{PROPOSITION 3.}{For any $G$ and for any $i,j=\1n$ such that
$j\ne i,\;$ $q\_{ii}>q\_{ij}$.
}

This property has a natural interpretation, namely, each vertex is more
``accessible" from itself than from any other vertex.

\statement{PROPOSITION~4~(triangle inequality for proximities).}%
{For any $G$ and for any $i,j,k=\1n,\;$ $q\_{ij}+q\_{ik}-q\_{jk}\le
q\_{ii}.$ If$,$ in addition$,$ $i\ne j$ and $i\ne k,$ then
$q\_{ij}+q\_{ik}-q\_{jk}< q\_{ii}.$
}

Consider the index
\begin{equation}
d\_{ij}=q\_{ii}+q\_{jj}- q\_{ij}-q\_{ji}=
        q\_{ii}+q\_{jj}-2q\_{ij},\quad i,j=\1n.
\label{metr}
\end{equation}

\statement{ASSERTION 1.}{$d(i,j)=d_{ij},\:$ $i,j=\1n,$ is a distance
function for multigraph vertices$,$ i.e.$,$ it complies with the axioms
of metric.
}

This assertion is easily proved using the above propositions (this is
left to the reader). The triangle inequality for proximities turns out
to be equivalent to the ordinary triangle inequality for metric
$d_{ij}$, which justifies the name of the former inequality. In
contrast to the standard graph distance, this metric considers all
connections in a graph.

\statement{PROPOSITION 5.}{For any $G$ and for any
$i,j=\1n,\;\:q\_{ij}=0$ iff there exist no paths between $i$ and~$j.$
}\hspace{-.3ex}

\secondstatement{COROLLARY.}
{{\rm (1)} Matrix $Q$ is reducible to a block-diagonal form, where all
block entries are strictly positive and all off-block entries are zeros.
$Q$ is strictly positive iff multigraph $G$ is connected.
\newline{\rm (2)} For any $i,j,k \in V(G),$ if $q\_{ij}>0$ and
$q\_{jk}>0,$ then $q\_{ik}>0$.
}

\secondstatement{PROPOSITION 6.}{For any $G$ and for any $i,k,t=\1n,$
\newline{\rm (1)} if there exists a path in $G$ from $i$ to $k,$
$t\ne k,$ and every path from $i$ to $t$ includes $k,$ then
$q\_{ik}>q\_{it}.$
\newline{\rm (2)} if $q\_{ik}>q\_{it}$ and $i\ne k,$ then there
exists a path from $i$ to $k$, such that the difference
$(q\_{jk}-q\_{jt})$ strictly increases as $j$ progresses from $i$
to $k$ along the path.
}

\secondstatement{PROPOSITION 7.}{Suppose that some edge weight
$\e_{kt}^p$ in $G$ increases by $\D\e\_{kt}>0$ or an extra edge between
$k$ and $t$ with a strictly positive weight $\D\e\_{kt}$ is added to
$G.$ Let $G'$ be the new graph and $W'=W(G'),$ $Q'=Q(G').$
Then
\newline
\indent{\rm (1)} $\D Q=hR,$ where $\D Q=Q'-Q,$
$h=\dfrac{\D\e\_{kt}}{\D\e\_{kt}(q\_{kk}+q\_{tt}-2q\_{kt})+1}
=(d\_{kt}+1/\D\e\_{kt})^{-1},$ and $R=(r\_{ij})$ is the $n\times n$
matrix with entries $r\_{ij}=(q\_{ik}-q\_{it})(q\_{jt}-q\_{jk});$
\newline
\indent{\rm (2) (}this item and the following three are corollaries
from item~{\rm (1))} all rows and all columns of $\D Q$ are
proportional$,$ i.e.$,$ ${\rm rank}\D Q=1;$
\newline
\indent{\rm (3)} if
$q\_{ik}>q\_{it},$ then $\D q\_{ij}>0$ iff
$q\_{jt}>q\_{jk},$ and $\D q\_{ij}<0$ iff
$q\_{jk}>q\_{jt};$
\newline
\indent{\rm (4)} the signs of all increments $\D q\_{ij}$ do not depend
on the absolute value of $\D\e\_{kt},$ and the absolute values of
nonzero $\D q\_{ij}$ strictly increase in $\D\e\_{kt};$
\newline
\indent{\rm (5)} for any $i,j\in V(G),\;$ $\D d\_{ij}=-{1\over 4}
(d\_{ik}-d\_{it}+d\_{jt}-d\_{jk})^2(d\_{kt}+1/\D\e\_{kt})^{-1},$ and
therefore $d'_{ij}\le d\_{ij}.$
}

According to item~(3), if the direct connection between $k$ and $t$
intensifies, then the relative accessibility of $j$ from $i$ increases
if and only if $i$ and $j$ initially were ``more strongly connected"
with {\em different\/} vertices of the pair $(k,t)$. Otherwise, it can
be said that the connections in the multigraph intensify outside of
most paths from $i$ to $j$, thus the relative accessibility of $j$ from
$i$ decreases.

Propositions~8 and~9 are corollaries of Proposition~7.

\statement{PROPOSITION 8.}{Suppose that some edge weight $\e_{kt}^p$ in
$G$ increases or an extra edge between $k$ and $t$ with a positive
weight is added to $G.$ Then
\newline
\indent{\rm (1)} $q\_{kt}$ increases$;$
\newline
\indent{\rm (2)} for any $i=\1n,$ if there exists a path from $i$ to
$k$ and every path from $i$ to $t$ includes $k,$ then
$\D q\_{it}>\D q\_{ik};$
\newline
\indent{\rm (3)} for any $i\_1,i\_2=\1n,$ if both $i\_1$ and $i\_2$ can
be substituted for $i$ in the hypothesis of item~{\rm(2),} then
$q\_{i\_1i\_2}$ decreases$;$
\newline
\indent{\rm (4)} for any $i=\1n,$ if $q\_{ik}=q\_{it},$ then $q\_{ij}$
do not alter for all $j=\1n$.
}

By item~(3), the relative accessibility between a pair of vertices
decreases when some ``extraneous" connections appear or intensify in
$G$.

Let $D$ be a subset of vertex set $V(G)$. We say that $D$ is a {\it
macrovertex in $G$} if for all $i\in D$, $j\in D$, and $k\notin D,\;\:
\ell\_{ik}=\ell\_{jk}$.

The following property is among to the most interesting ones. It
provides a sufficient condition for the equality and stability of
relative forest accessibilities.

\statement{PROPOSITION~9~(macrovertex independence).}{ Suppose that
$D$ is a macrovertex in $G$ and $i\in D,$ $j\in D,$ $k\notin D.$ Then
\newline
\indent{\rm (1)} $q\_{ik}=q\_{jk};$
\newline
\indent{\rm (2)} $q\_{ik}$ does not alter when any new edges appear or
the weights of any existing edges change inside $D.$
}

Now we shall obtain an alternative topological interpretation of the
matrix $Q$ of relative forest accessibilities (the first interpretation
is provided by Theorem~3). It will be demonstrated that $q\_{ij}$ are
related to the weights of routes of various lengths between $i$ and $j$
in $G$. To be more precise, introduce the notion of route with drains.

A {\it route with drains\/} (RWD) is an alternating sequence of
multigraph vertices and edges with the following features:

(1) the sequence starts and ends with vertices;

(2) the edge located between two {\it different\/} vertices in the
sequence is incident to them. If the same vertex stands in the sequence
before and after an edge, it is only required that it be incident to
this edge, the second incident vertex being arbitrary. Such an edge is
called a {\it drain}.

Routes with drains result from the usual routes by adding any number of
one-edge offshoots (drains), which may, in particular, follow ``forward"
and ``backward" along the original route.

The total number of edges in the sequence is called the {\it length\/}
of the route with drains. Set, by definition, te fact that for any
vertex $i$ there exist one route of length $0$ from $i$ to $i$ with $0$
drains and no other routes of length $0$.

The {\em weight\/} of a route with drains is defined as the product of
the weights of all its edges (if an edge enters a route
with drains $k$ times, its weight is taken with exponent $k$). For
any $i=\1n$, the weight of the $0$-length RWD from $i$ to $i$ is
set to be $1$.

Let $a^*=\maxl_{i,j\in V(G)}a\_{ij}$ be the maximal number of
multiple edges incident to any pair of vertices in $G$.

\statement{PROPOSITION 10.}{For any weighted multigraph $G$ with all
weights of edges from the interval $(0,\big(2a^*(n-1)\big)^{-1})$ and
for any $i,j=\1n,$
\[
q\_{ij}=\suml_{t=0}^\infty (U_{ij}^{(t)}-P_{ij}^{(t)}),
\]
where $U_{ij}^{(t)}$ and $P_{ij}^{(t)}$ are the total weights of all
routes of length $t$ with even and odd number of drains between
vertices $i$ and $j$ in $G,$ respectively.
}

One more interpretation of $Q$ can be obtained with the help of the
Cayley--Hamilton theorem \cite{GoDrRo}.

Instead of $Q$, one can use the matrices $Q_{\alpha}=(I+\alpha
L)^{-1}$, $\alpha>0$, which have the same properties as $Q$ except for
Proposition~10, where the factor $\alpha$ appears. The parameter
$\alpha>0$ specifies the proportions of accounting for long connections
between vertices of $G$ versus short ones.

\arcsect{4}{ACCESSIBILITY AND DERIVATIVE STRUCTURAL INDICES}

The foregoing properties of the relative forest accessibility
demonstrate that it is an appropriate index of proximity (connectivity,
accessibility) of graph (multigraph) vertices. A distinctive feature of
this index is its normalization: the sum of the accessibilities of all
vertices from a given one and the sum of the accessibilities of a given
vertex from all vertices of a multigraph are equal to unity. Therefore,
each $i$th row of the matrix $Q$ can be treated as a probability
distribution (or shares of a certain resource) somehow related to the
vertex $i$. In which cases is such a normalization necessary? Consider
two examples.

Suppose that the members of a group collect information from the
environment and exchange it with each other, the intensity of the
exchange being specified for each pair. Every participant transmits not
only the information collected on her own, but also that received from
the others. It is required to ascertain which fractions of the
cumulative information received by the $i$th participant were initially
collected by each member of the group. In this example, information can
be replaced with, for example, influence or material resources. The
principal feature is the distribution of some resource related to a
certain vertex, over all vertices.

The second example is a variant of a children's ``ring" game in which
the ring may successively be passed many times, and this is done
secretly, not before the players' eyes. If the pairwise transfer
probabilities are specified for all players along with the temporal
parameters of this random process (which is a Markov process in the
simplest case), one can take an interest in the ring's location
probabilities at every moment, provided that its starting location was
at vertex (player) $i$. The main feature of this example is the
presence of probability distributions related to each vertex.

In the above examples, if an adequate mathematical model is stated, the
result is precise, not heuristic, and there is no need to select it
being guided by ``good properties," such as those given in the previous
section. It turns out, however, that for both examples there are
rather natural models (we intend to describe them and compare them with
other models, e.g., \cite{Frie,Yama}, in our next paper), which lead to
relative forest accessibilities. This means, in turn, that even when
there is no detailed model, only the intensities (or probabilities) of
pairwise interactions being known, the relative forest accessibilities
provide a comprehensible first approximation for the required values.

Now we turn to derivative structural indices. The value
\[
1-\suml_{j\ne i}q\_{ij}=q\_{ii}
\]
can serve to measure the {\em solitariness\/} of the $i$th member of a
group. Now a number of other indices can be constructed in the usual
fashion. Specifically, the mean solitariness over a group,
\[
\rho=\frac{1}{n}\suml_{i=1}^n q\_{ii},
\]
indicates the extent of its {\it dissociation}. The empirical variance
of the solitariness evaluates the {\it heterogeneity\/} of a group. The
ratio of $q\_{ii}$ to $\rho$ (or their difference) measures the
{\em provinciality\/} of the $i$th member of a group.
Equation~(\ref{metr}) introduces a specific distance between the
members of a group (Assertion~1 in Section~3%
). The properties of all these indices are determined by those of the
relative forest accessibilities studied above.

Notice, in conclusion, that there exists a certain relation between the
problem of centrality (respectively, provinciality) evaluation and the
problem of estimating the strength of players from incomplete
tournaments. In the latter case, an ``object--object" matrix is
processed as well, but its entries express the results of {\em paired
comparisons\/} (e.g., games or comparative preferences) rather than
personal choices within a group. The problem of scoring from paired
comparisons has been investigated a little bit better (but also
insufficiently). It is worth noting, for example, that the work
\cite{Kat} was accepted as relevant in the literature on paired
comparisons, though it was concerned with sociometric data. And
conversely, sensitive scoring methods for preference aggregation can
be considered with reference to sociometric data. A review of these
methods can be found in \cite{ChSh}.

\arcpril  

\PLE{1} Lemma~1 is reducible to Lemma~2, since for every multigraph
$G$, the corresponding multidigraph $\G$ can be introduced by replacing
every edge of $G$ with a pair of opposite arcs with the same weight
each. The matrix $W$ (and thus $Q$) is the same for $G$ and $\G$, so
the desired statements of Lemma~1 follow from the existence of a
natural one-to-one correspondence between all spanning rooted forests
in $G$ and all spanning diverging forests in $\G$.

{\bf Proposition 1} follows from the symmetry of $W$.

\PPR{2} Item~(1) follows from Theorem~3 and the positiveness of edge
weights.

Item~(2) immediately follows from the fact that $W=Q^{-1}$ satisfies the
same condition \cite{GoDrRo,Sha}. Another easy proof is provided by
Theorem~3 and the fact that for any $i\_1,i\_2,j\in V(G)$, $i\_1\ne
i\_2\;\Rightarrow \FF^{i\_1j}\cap \FF^{i\_2j}=\varnothing$ and
$\cupl^n_{i=1}\FF^{ij}=\FF$.

Item~(3) follows from item~(2) and Proposition~1.

\PPR{3} Note that for any $i,j=\1n$ such that $j\ne i$ and for any
$H\in\FF$, if $H\in\FF^{ij}$ then $H\in\FF^{ii}.$ Therefore, $\FF^{ij}
\subseteq \FF^{ii}$. Let $F\_0$ be a subgraph of $G$ such that
$V(F\_0)=V(G)$ and $E(F\_0)=\varnothing$. Then
$F\_0\in\setminus{\FF^{ii}}{\FF^{ij}}$ and $\e(F\_0)=1$, i.e.,
$\FF^{ij}\subset\FF^{ii}$ and $\e(\FF^{ij})<\e(\FF^{ii})$. By
Theorem~3, $q\_{ii}>q\_{ij}$.

\PPR{4} If $i=j$ or $i=k$ then, obviously,
\[
q\_{ij}+q\_{ik}-q\_{jk}=q\_{ii}.
\]

Assume that $i\ne j$ and $i\ne k.$ In the same way as in the proof of
Proposition~3, we have
\[
\FF^{ij}\cup\FF^{ik}\subset\FF^{ii},
\]
and hence
\begin{equation}
 \e(\FF^{ij}\cup\FF^{ik})
=\e(\FF^{ij})+\e(\FF^{ik})-\e(\FF^{ij}\cap\FF^{ik})<\e(\FF^{ii}).
\label{simra}
\end{equation}

Define $\FF^{ijk}$ as $\FF^{ij}\cap\FF^{ik}$. Observe that $\FF^{ijk}$
differs from $\FF^{jik}=\FF^{ji}\cap\FF^{jk}$ only by the roots in the
trees containing $i$, $j$, and $k$ simultaneously. Therefore,
\begin{equation}
\e(\FF^{ij}\cap \FF^{ik})=\e(\FF^{ijk})=\e(\FF^{jik})\le\e(\FF^{jk}).
\label{treuh}
\end{equation}

Inequalities (\ref{simra}) and (\ref{treuh}) imply
\[
\e(\FF^{ij})+\e(\FF^{ik})-\e(\FF^{jk})<\e(\FF^{ii}),
\]
and, by Theorem 3,
\[
q\_{ij}+q\_{ik}-q\_{jk}<q\_{ii}.
\]

{\bf Proposition 5} follows directly from Theorem~3.

\PPR{6} Item~(1). Note that $H\in\FF^{it}$
implies $H\in\FF^{ik}$. On the other hand,
$\setminus{\FF^{ik}}{\FF^{it}}\ne\varnothing$ and
$\e(\setminus{\FF^{ik}}{\FF^{it}})>0$. Hence by Theorem~3,
$q\_{ik}>q\_{it}$.

Item~(2). By virtue of Eq.~(\ref{DefQ}),
\begin{equation}
(I+L)Q=I.
\label{inver}
\end{equation}

Rewrite (\ref{inver}) componentwise for entries $ik$ and $it$ of the
matrix $(I+L)Q$. Using Eqs.~(\ref{lij}) and (\ref{lii}), the notation
$\e\_{ij}=-\l\_{ij}$, and $i\ne t$, which follows from Proposition~3,
we get
\[
q\_{ik}=\suml_{j\ne i}\e\_{ij}(q\_{jk}-q\_{ik}),
\]
\[
q\_{it}=\suml_{j\ne i}\e\_{ij}(q\_{jt}-q\_{it}),
\]
\[
q\_{ik}-q\_{it}=\suml_{j\ne i}\e\_{ij}[(q\_{jk}-q\_{jt})-
(q\_{ik}-q\_{it})].
\]

Then, since $q\_{ik}-q\_{it}>0$, there exists $j\ne i$ such that
$\e\_{ij}\ne 0$ (and thus $(ij)\in E(G)$) and
$q\_{jk}-q\_{jt}>q\_{ik}-q\_{it}$ (recall that the case $\e\_{ij}<0$ is
excluded).

Applying this argument to vertex $j$ instead of $i$, and so forth, and
taking into account that no vertex in the path thereby constituted may
coincide with any previous one and that $i\ne k$, we finally obtain $k$
as the terminal vertex of this path, as desired.

\PPR{7} Let $\D W=W'-W$. Note that $\D W =XY$, where $X=(x\_{i1})$,
$i=\1n$, is the column vector with entries $x\_{k1}=1$, $x\_{t1}=-1$,
and $x\_{i1}=0$ for all $i\ne k$, $i\ne t$; $Y=(y\_{1j})$, $j=\1n$, is
the row vector with entries $y\_{1k}=\D\e\_{kt}$,
$y\_{1t}=-\D\e\_{kt}$, and $y\_{1j}=0$ for all $j\ne k$, $j\ne t$.
According to \cite[Sec.~0.7.4]{HoJo},
\[
Q'=Q-\frac{1}{1+YQX}QXYQ.
\]
It is straightforward to verify that $(-\tfrac{1}{1+YQX})=-h/\D\e\_{kt}$
and $QXYQ=-\D\e\_{kt}R$, and thereby item~(1) is proved. Items~(2)
through (5) follow from item~(1) and the nonnegativity of $d\_{kt}$
(see Proposition~3 or Assertion~1).

\PPR{8} Item~(1). By Proposition~3, $q\_{kk}>q\_{kt}$ and
$q\_{tt}>q\_{tk}$, and hence item~(3) of Proposition~7 implies $\D
q\_{kt}>0$.

Item~(2). Setting $Q'=Q(G')$, by item~(1) of Proposition~7 we have
\begin{eqnarray*}
\D q\_{it}-\D q\_{ik}
\I&=&\I h(q\_{ik}-q\_{it})(q\_{tt}-q\_{tk})
  -h(q\_{ik}-q\_{it})(q\_{kt}-q\_{kk})
\\
\I&=&\I h(q\_{ik}-q\_{it})(q\_{kk}+q\_{tt}-q\_{tk}-q\_{kt})
=h(q\_{ik}-q\_{it})d\_{kt}.
\end{eqnarray*}
Now the desired inequality follows from item~(1) of Proposition~6
together with Assertion~1.

Item~(3). By item~(1) of Proposition~6, $q\_{i_1k}>q\_{i_1t}$ and
$q\_{i_2k}>q\_{i_2t}$, and by item~(3) of Proposition~7,
$\D q\_{i_1i_2}<0$.

Item~(4). By item~(1) of Proposition~7 we have
$\D q\_{ij}=h(q\_{ik}-q\_{it})(q\_{jt}-q\_{jk})=0.$

\PPR{9} Consider the graph $G'$ on the vertex set $V(G)$ such that\\
\indent(1) $(ij)\in E(G')$ iff $i\ne j$ and $\l\_{ij}\ne 0$, and\\
\indent(2) for every edge $(ij)\in E(G'),\;$ $\e'_{ij}=-\l\_{ij}$.

Let $Q'=Q(G')=(q'_{ij})$. Obviously, $D$ is a macrovertex in $G'$ as
well as in $G$. Let $S=\setminus{V(G)}{D}$. First, we prove
Proposition~9 for $G'$. Consider the graph $G''$ resulting from $G'$ by
deleting all edges inside $D$. Let $Q''=Q(G'')=(q''_{ij})$. All
vertices of $D$ are symmetric in $G''$; therefore, $q''_{ik}=q''_{jk}$
for any $i,j\in D$, $k\in S$. Then using item~(4) of Proposition~8, by
means of induction we get $q'_{ik}=q''_{ik}=q''_{jk}=q'_{jk}$, for all
$i,j\in D$ and $k\in S$. This proves Proposition~9, since $Q'=Q$.

\PPR{10} Expand $Q=(I-(-L))^{-1}$ as the sum of an infinitely
decreasing geometric progression using the notation $M=(m\_{ij})=-L$:
\begin{equation}
Q=(I-M)^{-1}=I+M+M^2+\ldots\quad.
\label{expand}
\end{equation}

This expansion is valid if and only if
\begin{equation}
\vert\la\_1\vert<1,
\label{lamb}
\end{equation}
where $\vert\la\_1\vert$ is the spectral radius of $M=-L$
\cite[Corollary 5.6.16]{HoJo}.

Consider the upper bound of $\vert\lambda\vert\_{\max}$ provided by
the Ger\v{s}gorin theorem (see \cite{HoJo}):
\begin{equation}
\vert\la\_1\vert\le\maxl_{1\le i\le n}\suml_{j=1}^n
\vert \l\_{ij}\vert.
\label{Gers}
\end{equation}
Let $\e\_{\max}=\maxl_{1\le i\ne j\le n}\e\_{ij}$, where
$\e\_{ij}=\suml_{p=1}^{a\_{ij}}\e_{ij}^p=-\l\_{ij}$. Then by
Eqs.~(\ref{lij}) and (\ref{lii}),
\begin{equation}
    \maxl_{1\le i\le n}\suml_{j=1}^n \vert\l\_{ij}\vert
  =2\maxl_{1\le i\le n}\suml_{j\ne i}\vert\l\_{ij}\vert
\le2\maxl_{1\le i\le n}\suml_{j\ne i}a^*\e\_{\max}
=2a^*(n-1)\e\_{\max}.
\end{equation}
Consequently, the fulfillment of (\ref{lamb}), and therefore of
(\ref{expand}) is assured by
\[
\e\_{\max}<\big(2a^*(n-1)\big)^{-1}.
\]

By virtue of (\ref{expand}), it suffices to prove that
\begin{equation}
m_{ij}^{(k)}=U_{ij}^{(k)}-P_{ij}^{(k)},\quad i,j=\1n,\quad
k=0,1,2,\ldots,
\label{route}
\end{equation}
where $m_{ij}^{(k)},\;i,j=\1n$, are the entries of $M^k$.

Let us apply induction on the length $k$ of the roots with drains
between $i$ and $j$. The proof can be used with no change for the case
of digraphs, because it does not use the symmetry of $M$.

$1^0$. $k=0$. Equation~(\ref{route}) is valid because $M^0=I$ and by
the definition of root with drains, for every $i,j=\1n,\;j\ne i,\;\;$
$U_{ii}^{(0)}=1$ and $P_{ii}^{(0)}=P_{ij}^{(0)}=U_{ij}^{(0)}=0$ hold.

$2^0$. Let (\ref{route}) be valid for $k=v$. Prove it for $k=v+1$.
Consider an arbitrary route $\mu$ of length $v+1$ with $g$ drains
between vertices $i$ and $j$. Let $t$ be the next to last vertex of
$\mu$. If $t\ne j$, then $\mu$ is representable as the combination of a
route of length $v$ with $g$ drains from $i$ to $t$ and an edge $(tj)$.
Otherwise, $t=j$, and $\mu$ can be considered as the combination of a
route between $i$ and $j$ with $g-1$ drains and an edge incident with
$j$ (this edge is the $g$th drain). Therefore,
\[
U^{(v+1)}_{ij}
=\suml_{t\ne j}U_{it}^{(v)}m\_{tj}
+\suml_{t\ne j}P_{ij}^{(v)}m\_{jt},
\]
\[
P^{(v+1)}_{ij}
=\suml_{t\ne j}P_{it}^{(v)}m\_{tj}
+\suml_{t\ne j}U_{ij}^{(v)}m\_{jt}.
\]

Then
\[
U^{(v+1)}_{ij}-P^{(v+1)}_{ij}
=\suml_{t\ne j} U_{it}^{(v)}m\_{tj}
+\suml_{t\ne j} P_{ij}^{(v)}m\_{jt}
-\suml_{t\ne j} P_{it}^{(v)}m\_{tj}
-\suml_{t\ne j} U_{ij}^{(v)}m\_{jt}
\]
\[
=\suml_{t\ne j} (U_{it}^{(v)}-P_{it}^{(v)})m\_{tj}
-\suml_{t\ne j} (U_{ij}^{(v)}-P_{ij}^{(v)})m\_{jt}\eql^{\langle1\rangle}
 \suml_{t\ne j} m_{it}^{(v)}m\_{tj}-m_{ij}^{(v)}
 \suml_{t\ne j} m\_{jt}
\]
\[
 \eql^{\langle2\rangle}
 \suml_{t\ne j} m_{it}^{(v)}m\_{tj}+m_{ij}^{(v)} m\_{jj}
=\suml_{t=1}^n m_{it}^{(v)}m\_{tj}
=m_{ij}^{(v+1)},
\]
where transition $\langle1\rangle$ is carried out by the induction
hypothesis, and $\langle2\rangle$ uses the equality
$m\_{jj}=-\suml_{t\ne j} m\_{jt}$ which follows from Eq.~(\ref{lii})
using $M=-L$. Proposition~10 is proved.


\begin{thebibliography}{99}

\bibitem{Pan} V.~I.~Paniotto, ``Analysis of the structure of
interpersonal relations," in: {\it Mathematical Methods of the Analysis
and Interpretation of Sociological Data} [in Russian], Nauka, Moscow
(1989), pp.~121--162.
\vs
\bibitem{HeGl} T.~H{\o}vik and N.~P.~Gleditsch, ``Structural parameters
of graphs: A theoretical investigation," in: H.~M.~Blalock,
A.~Aganbegian, F.~M.~Borodkin, R.~Boudon, and V.~Capecchi (eds.)
{\it Quantitative Sociology. International Perspectives on Mathematical
and Statistical Modeling}, Academic Press, New York (1975),
pp.~203--223.
\vs
\bibitem{ChShAtl} P.~Yu.~Chebotarev and E.~Shamis, ``On the proximity
measure for graph vertices provided by the inverse Laplacian
characteristic matrix," in: {\it 5th Conference of the International
Linear Algebra Society}, Georgia State University, Atlanta (1995),
pp.~30--31.
\vs
\bibitem{Che89} P.~Yu.~Chebotarev, ``Generalization of the row sum
method for incomplete paired comparisons," {\it Automat. Remote
Control}, {\bf 50}, No.~8, Part~2, 1103--1113 (1989).
\vs
\bibitem{ShaAMS} E.~Shamis, ``Counting spanning converging forests,"
{\it Abstracts of Papers Presented to the American Mathematical
Society}, {\bf 15}, 412--413 (1994).
\vs
\bibitem{Sha} E.~Shamis, ``Graph-theoretic interpretation of the
generalized row sum method," {\it Math. Soc. Sci.}, {\bf 27}, 321--333
(1994).
\vs
\bibitem{Merr} R.~Merris, ``Doubly stochastic graph matrices,"
{\it Publikacije Elektrotehn. Fakulteta, Univerzitet U Beogradu},
{\bf 8}, 6--13 (1997).
\vs
\bibitem{GoDrRo} V.~E.~Golender, V.~V.~Drboglav, and A.~B.~Rosenblit,
``Graph potentials method and its application for chemical information
processing," {\it J. Chem. Inf. Comput. Sci.}, {\bf 21}, 196--204
(1981).
\vs
\bibitem{Frie} N.~E.~Friedkin, ``Theoretical foundations of centrality
measures," {\it Amer. J. Sociology}, {\bf 96}, 1478--1504 (1991).
\vs
\bibitem{Moo} J.~W.~Moon, {\it Counting Labelled Trees}, Canad. Math.
Congress, Montreal (1970).
\vs
\bibitem{Tut} W.~T.~Tutte, {\it Graph Theory}, Addison-Wesley,
Reading, Massachusetts (1984).
\vs
\bibitem{HaPa} F.~Harary and E.~M.~Palmer, {\it Graphical Enumeration},
Academic Press, New York (1973).
\vs
\bibitem{Cha} S.~Chaiken, ``A combinatorial proof of the all minors
matrix tree theorem," {\it SIAM J. Alg. Disc. Meth.}, {\bf 3}, 319--329
(1982).
\vs
\bibitem{MoonDM} J.~W.~Moon, ``Some determinant expansions and the
matrix-tree theorem," {\it Discrete Math.}, {\bf 124}, 163--171 (1994).
\vs
\bibitem{Kel} A.~K.~Kelmans, ``On the properties of the characteristic
polynomial of a graph," in: {\it Cybernetics Serves Communism} [in
Russian], Vol.~4, Energiya, Moscow--Leningrad (1967), pp.~27--41.
\vs
\bibitem{KeCh} A.~K.~Kelmans and V.~M.~Chelnokov, ``A certain
polynomial of a graph and graphs with an extremal number of trees,"
{\it J. Comb. Theory}, Ser.~B, {\bf 16}, 197--214 (1974).
\vs
\bibitem{FiSe} M.~Fiedler and J.~Sedl\'{a}\v{c}ek, ``O $W$-bas\'{\i}ch
orientovan\'{y}ch graf\uu," {\it \v{C}asopis P\v{e}st. Mat.}, {\bf 83},
214--225 (1958).
\vs
\bibitem{MaOl} J.~S.~Maybee, D.~D.~Olesky, P.~van~den~Driessche, and
G.~Wiener, ``Matrices, digraphs and determinants," {\it SIAM J.
Matrix Anal. Appl.}, {\bf 10}, 500--519 (1989).
\vs
\bibitem{Yama} K.~Yamaguchi, ``The flow of information through
social networks: diagonal-free measures of inefficiency and the
structural determinants of inefficiency," {\it Social Networks},
{\bf 16}, 57--86 (1994).
\vs
\bibitem{Kat} L.~Katz, ``A new status index derived from sociometric
analysis," {\it Psychometrika}, {\bf 18}, 39--43 (1953).
\vs
\bibitem{ChSh} P.~Yu.~Chebotarev and E.~Shamis, ``Preference fusion
when the number of alternatives exceeds two: Indirect scoring
procedures," in: {\it Proceedings of Workshop on Foundations of
Information / Decision Fusion}, ed. by N.~S.~V.~Rao, V.~Protopopescu,
J.~Bernen, and G.~Seetharaman, Acadiana, Lafayette (1996), pp.~20--32.
\vs
\bibitem{HoJo} R.~A.~Horn and C.~R.~Johnson, {\it Matrix Analysis},
Cambridge University Press, Cambridge (1986).

\end{thebibliography}
\end{document}